\documentclass[12pt]{article}%
\usepackage[applemac]{inputenc}
\usepackage{amsmath,amssymb}
\usepackage{amsfonts}
\usepackage{amsmath,amssymb}
\usepackage[applemac]{inputenc}
\usepackage{color}
\usepackage{color, colortbl}
\usepackage{amsmath}
\usepackage{amssymb}
\usepackage{bbm}
\usepackage{graphicx}
\usepackage{tikz}
\usepackage{geometry}
\usetikzlibrary{arrows}
\usepackage{amsfonts}
\usepackage{amsmath,amssymb}
\usepackage[applemac]{inputenc}
\usepackage{color}
\usepackage{amsmath}
\usepackage{amssymb}
\usepackage{graphicx}
\usepackage{tikz}
\newtheorem{theorem}{Theorem}[section]

\newtheorem{definition}[theorem]{Definition}
\newtheorem{notation}[theorem]{Notation}
\newtheorem{example}[theorem]{Example}

\newtheorem{problem}[theorem]{Problem}
\newtheorem{remark}[theorem]{Remark}

\usetikzlibrary{arrows,shapes,automata,petri}

\newcommand{\dproof}{\noindent {Proof.} \quad}
\newcommand{\fproof}{\hfill $\square$ \bigskip}

\numberwithin{equation}{section}

\def\RB{\mathbb{R}}

\def\FB{\mathbb{F}}

\definecolor{LightCyan}{rgb}{0.88,1,1}

\def\MM{{\mathbb{M}}}

\def\EE{{\mathbb{ E}}}

\def\1B{\text{1\!\!I}}

\def\l{\langle}
\def\r{\rangle}
\def\<{\langle}
\def\>{\rangle}
\def\Z{\mathcal{Z}}
\def\F{\mathcal{F}}
\def\V{\mathcal{V}}

\def\H{\mathcal{H}}
\def\M{\mathcal{M}}
\def\T{\mathcal{T}}
\def\S{\mathcal{S}}
\def\L{\mathcal{L}}
\def\V{\mathcal{V}}
\def\G{\mathcal{G}}
\def\C{\mathcal{C}}

\begin{document}

\title{Impulse control of conditional McKean-Vlasov jump diffusions}
\author{Nacira Agram$^{1},$ Giulia Pucci$^{1}$ \& Bernt \O ksendal$^{2}$}
\date{\today}
\maketitle

\footnotetext[1]{Department of Mathematics, KTH Royal Institute of Technology 100 44, Stockholm, Sweden. \newline
Email: nacira@kth.se, pucci@kth.se. Work supported by the Swedish Research Council grant (2020-04697).}

\footnotetext[2]{%
Department of Mathematics, University of Oslo, Norway. 
Email: oksendal@math.uio.no.}

\begin{abstract}
This paper establishes a verification theorem for impulse control problems involving conditional McKean-Vlasov jump diffusions. 
We obtain a Markovian system by combining the state equation of the problem with the stochastic Fokker-Planck equation for the conditional probability law of the state. 
We derive sufficient variational inequalities for a function to be the value function of the impulse control problem, and for an impulse control to be the optimal control. 
We illustrate our results by applying them to the study of an optimal stream of dividends under transaction costs. We obtain the solution explicitly by finding a function and an associated impulse control which satisfy the verification theorem.
\end{abstract}



\section{Introduction}
Consider a filtered probability space $(\Omega,\F,P, \FB=\{\F\}_{t\geq 0})$  on which we are given a $d$-dimensional Brownian motion $B=(B_1, B_2, \ldots , B_d)$, a $k$-dimensional compensated Poisson random measure $\widetilde{N}(dt,d\zeta)$ such that 
$$ \widetilde{N}(dt,d\zeta)=N(dt,d\zeta)-\nu(d\zeta)dt, $$
where $N(dt,d\zeta)$ is a Poisson random measure and $\nu(d\zeta)$ the L\'evy measure of $N$,  and a random variable $Z \in L^2(P)$ that is independent of $\FB$. We denote by $L^2(P)$ the set of all the $d-$dimensional $\mathcal{F}$-measurable random variables $X$ such that $\EE[X^2]<\infty,$ 
where $\mathbb{E}$ denotes the expectation with respect to $P$. We consider the state process $X(t) \in \RB^d$ given as the solution of the following \textit{conditional McKean-Vlasov jump equation}
\small
\begin{align}
X(t) &=Z+\int_0^t\alpha(s,X(s),\mu_s)dt+\beta(s,X(s),\mu_{s})dB(s)\nonumber\\
&+\int_0^t\int_{
\mathbb{R}^d
}\gamma(s,X(s^-),\mu_{s^{-}},\zeta)\widetilde{N}(ds,d\zeta),\label{MV}
\end{align}
where we denote by 
$\mu_t=\mathcal{L}(X(t) | \mathcal{F}_t^{(1)})$
 the conditional probability distribution of $X(t)$ given the filtration $\mathcal{F}_t^{(1)}$generated by the first component $B_1(u);u\leq t$ of the Brownian motion up to time $t$.
Loosely speaking, the equation above models a McKean-Vlasov dynamics which is subject to what is called a "common noise" coming from the Brownian motion $B_1(t)$, which is observed and is influencing the dynamics of the system.\\
So defined, $\mu_t$ is a Borel probability measure on $\RB^d$ 
for all $t \in [0,T],\; \omega \in \Omega$. In particular, $\mu_t \in \MM_0$, with $\MM_0$ the set of deterministic Radon measures i.e. Borel measures finite on compact sets, outer regular on all Borel sets and inner regular on all open sets. Notice that all Borel probability measures on $\RB^d$ are Radon measures. From now on we will indicate with $\MM$ the set of random measures $\lambda(dx,\omega)$ which are Radon measures  with respect to $x$ for each $\omega$. We refer to \cite{F} for more information.\\
We suppose that 
$\alpha(t,x,\mu) \colon [0,T]\times\RB^d\times\MM \rightarrow \RB^d$, \; $\beta(t,x,\mu) \colon [0,T]\times\RB^d\times\MM \rightarrow \RB^{d\times m}$, \;$\gamma(t,x,\mu,\zeta) \colon [0,T]\times\RB^d\times\MM\times\RB^d \rightarrow \RB^{d\times k}$ are bounded processes and $\FB$-predictable for all $x,\mu,\zeta$ and they are also continuous with respect to $t$ and $x$ for all $\mu,\zeta$. \\
We can easily see that, under hypothesis of Lipschitz continuity and at most linear growth, there exists a unique solution for \eqref{MV} 
for all $t$ in $[0,T]$.

The purpose of this paper is to study impulse control problems for conditional McKean-Vlasov jump diffusions. In particular, we will define a performance criterion and then attempt to find a policy that maximizes performance within the admissible impulse strategies. Using a verification theorem approach, we establish a general form of quasi-variational inequalities and identify the sufficient conditions that lead to an optimal function. See precise formulation below. Standard impulse control problems can be solved by using the Dynkin formula. We refer to e.g. Bensoussan \& Lions \cite{BL} in the continuous case and to  \O ksendal and Sulem \cite{OS} in the setting of jump diffusions. \\ 
Impulse control problems naturally arise in many concrete applications, in particular when an operator, because of the intervention costs, decides to control the system by intervening only at a discrete set of times with a chosen intervention size: a sequence of stopping times $(\tau_1,\tau_2,\ldots,\tau_k, \ldots)$ is chosen to intervene and exercise the control. At each time $\tau_k$ of the player's $\textit{k}^{th}$ intervention, the player chooses an intervention of size $\zeta_k$. The impulse control consists of the sequence $\{ (\tau_k,\zeta_k) \}_{k\ge1}$. \\
Impulse control has sparked great interest in the financial field and beyond. See, for example, \cite{K} for portfolio theory applications, \cite{B} for energy markets, and \cite{CCTZ} for insurance. All of these works are based on quasi-variational inequalities and employ a verification approach.\\
Despite its adaptability to more realistic financial models, few papers have studied the case of mean field problems with impulse control. We refer to \cite{BCG} for a discussion of a more special type of impulse, where the only type of impulse is to add something to the system. This is a mean field game (MFG) where the mean-field (only the empirical mean) appears as an approximation of the many-player game. 
They
use the smooth fit principle (as used in the present work) to solve a specific MFG explicitly. \\
We refer also to \cite{CNS} for a MFG impulse control approach. Specifically, a problem of optimal harvesting in natural resource management is addressed. \\
A maximum principle for regime switching control problem for mean-field jump diffusions is studied by \cite{LLWH} but in that paper the problem considered is not really an impulse control problem because the intervention times are fixed in advance. \\
In our setting, we will not consider a MFG setup, as in the above mentioned works, we will only consider a decision maker who chooses the control to optimise a certain reward.  Moreover, the mean-field appears as a conditional probability distribution and to overcome the lack of the Markov property, we introduce the equation of the measure which is of stochastic Fokker-Planck type.\\
In \cite{DHH}, the authors can handle a non-Markovian dynamics. However, the impulse control is given in a particular compact form and only a given number of impulses are allowed. They use a Snell envelope approach and related reflected backward stochastic differential equations.\\
In the next section, we introduce some notations and present some preliminary results. As part of Section $3$, we state the optimal control problem and prove the verification theorem. In Section 4, we apply the previous results to solve an explicit problem of optimal dividend streams under transaction costs.

\section{Preliminaries}
 The process $X(t)$ given by \eqref{MV} is not in itself Markovian, so to be able to use the Dynkin formula, we extend the system to the process $Y$ defined by 
\begin{equation*}
Y(t)=(s+t,X(t),\mu_t); \quad t\geq 0;\quad Y(0)=(s,Z,\mu_0)=:y,
\end{equation*}
 for some arbitrary starting time $s\geq 0$, with state dynamics given by $X(t)$, conditional law of the state given by $\mu_t$ and with $X(0)=Z, \; \mu_0=\L(X(0))$. This system is Markovian, in virtue of the following Fokker-Planck equation for the conditional law 
 $\mu_t$, proved in \cite{AO}.

\begin{theorem}{(Conditional stochastic Fokker-Planck equation)} \label{FP1}\\
Let $X(t)$ be as in \eqref{MV} and let $\mu_t=\mu_t(dx,\omega)$ be the regular conditional distribution of $X(t)$ given $\mathcal{F}_t^{(1)}$. Then $\mu_t$ satisfies the following SPIDE (in the sense of distributions):
 
\begin{align} \label{FPmu}
d\mu _{t} =A_0^{*} \mu_t dt + A_1^{*}\mu_t dB_1(t);   \quad \mu_0=\mathcal{L}(X(0)),
\end{align}
where 
$A_0^{*}$ is the integro-differential operator
\begin{align*}
A_0^{*}\mu&= -\sum_{j=1}^d D_j [\alpha_j \mu] +\frac{1}{2}\sum_{n,j=1}^d D_{n,j}[(\beta \beta^{(T)})_{n,j} \mu] \nonumber\\
&+\sum_{\ell=1}^k  \int_{\RB}\Big\{\mu^{(\gamma^{(\ell)})}-\mu+\sum_{j=1}^d D_j[\gamma_j^{(\ell)}(s,\cdot,\zeta)\mu]\Big\} \nu_{\ell} \left( d\zeta \right), 
\end{align*}
and $A_1^{*}$ is the differential operator
\begin{align*}
A_1^{*}\mu= - \sum_{j=1}^d D_j[\beta_{1,j} \mu], 
\end{align*}
where $\beta^{(T)}$ denotes the transposed of the $d \times m$ - matrix $\beta=\big[\beta_{j,k}\big]_{1\leq j \leq d,1 \leq k \leq m}$ and $\gamma^{(\ell)}$ is column number $\ell$ of the matrix $\gamma$.

\end{theorem}

For notational simplicity,  we use $D_j, D_{n,j}$
to denote $\frac{\partial }{\partial x_j}$  and  $\frac{\partial^2}{\partial x_n \partial x_j}$  in the sense of distributions. \\
We have also used the following notation, taken from \cite{AO}.\\
For fixed $t,\mu,\zeta$ and $\ell=1,2,... k$, we write for simplicity $\gamma^{(\ell)}=\gamma^{(\ell)}(t,x,\mu,\zeta)$ for column number $\ell$ of the $d \times k$-matrix $\gamma$. Then $\nu_\ell$ represents the L\'evy measure of $N_\ell$ for all $\ell$. Note that for given $\mu \in \MM$ the map
$$g \mapsto \int_{\RB^d} g(x+\gamma^{(\ell)}) \mu(dx)$$
 is a bounded linear map on $C_0(\RB^d),$ which is defined to be the uniform closure of the space $C_c(\mathbb{R}^d)$ of continuous functions with compact support.
Therefore, since $\mathbb{M}$ is the dual of $\mathcal{C}_0(\mathbb{R}^d)$, there is a unique measure $\mu^{(\gamma^{(\ell)})} \in \MM$ such that 
\begin{align*}
\langle \mu^{(\gamma^{(\ell)})},g \rangle:=\int_{\RB^d} g(x) \mu^{(\gamma^{(\ell)})}(dx)=\int_{\RB^d}g(x+\gamma^{(\ell)}) \mu(dx), \text{ for all } g \in C_0(\RB^d),
\end{align*}
where $\langle\mu^{(\gamma^{(\ell)})},g \rangle$ denotes the action of the measure $\mu^{(\gamma^{(\ell)})}$ on $g$. We call $\mu^{(\gamma^{(\ell)})}$ the $\gamma^{(\ell)}$-shift of $\mu$. Note that $\mu^{(\gamma^{(\ell)})}$ is positive
 and absolutely continuous with respect to $\mu$.\\

\section{A General Formulation and a Verification Theorem}
As noted above, in virtue of the Fokker-Planck equation \eqref{FPmu} we can extend the system \eqref{MV} into a Markovian system by defining the following $[0,\infty) \times L^2(P) \times  \MM$ - valued process $Y(t):=(s+t, X(t),\mu_t)$ as follows:
\begin{align}
dY(t)&=F(Y(t))dt + G(Y(t))dB(t) + \int_{\RB^k}
H(Y(t^-),z)\widetilde{N}(dt,dz) \nonumber \\
&:=\left[\begin{array}{clcr}
dt \\dX(t)\\d\mu_t
\end{array} \right] =\left[ \begin{array}{c}
1\\ \alpha(Y(t)) \\A_0^{*}\mu_t
\end{array} \right] dt
+\left[\begin{array}{rc}
0_{1\times m} \\ \beta(Y(t)) \\ A_1^{*}\mu_t,0,0 ...,0
\end{array} \right]dB(t)\nonumber\\
&+ \int_{\RB^d} \left[ \begin{array}{rc}
0_{1\times k}\\ \gamma(Y(t^-),\zeta) \\0_{1\times k}
\end{array}\right] \widetilde{N}(dt,d\zeta),\quad s\leq t \leq T, \label{Y2}
\end{align}
where $X(t)$ and $\mu_t$ satisfy the equations \eqref{MV} and \eqref{FPmu}, respectively. 
Moreover, we have used the shorthand notation
\begin{align*}
\alpha(Y(t))&=\alpha(s+t,X(t),\mu(t))\\
\beta(Y(t))&= \beta(s+t,X(t),\mu(t))\\
\gamma(Y(t^-),\zeta)&=\gamma(s+t,X(t^-),\mu(t^-),\zeta).
\end{align*}
The process $Y(t)$ starts at $y=(s, Z,\mu)$. We shall denote by $\mu$ the initial probability distribution $\mathcal{L}(X(0))$ or the generic value of the conditional law $\mu_t :=\mathcal{L}(X(t) | \mathcal{F}_t^{(1)})$, when there is no ambiguity. Similarly, we use the following notation:
\begin{notation}
We use 
\begin{itemize}
\item
$x$ to denote a generic value of the point $X(t,\omega) \in \RB^d$, and
\item
$X$ to denote a generic value of the random variable $X(t, \cdot) \in L^2(P).$
\item
When the meaning is clear from the context we use $x$ in both situations.
\end{itemize}
\end{notation}

The concept of impulse control is simple and intuitive: at any time the agent can make an intervention $\zeta$ into the system. Due to the cost of each intervention the agent can intervene only at discrete times $\tau_1, \tau_2, \ldots $. The impulse problem is to find out at what times it is optimal to intervene and what is the corresponding optimal intervention sizes. We now proceed to formulate precisely our impulse control problem for conditional McKean-Vlasov jump diffusions.

Suppose that -- if there are no interventions --   the  $[0,\infty) \times L^2(P) \times  \MM$ - valued process
$Y(t)=(s+t,X(t),\mu_t)$ is the conditional McKean-Vlasov jump diffusion given by \eqref{Y2}.

Suppose that at any time $t$ and any state $y=(s,X,\mu)$ we are free to
intervene and give the state $X$ an impulse $\zeta\in\Z\subset\RB^d$,
where $\Z$ is a given set (the set of admissible impulse values).
Suppose the result of giving the state $X$ the impulse $\zeta$ is that the state jumps immediately from $X$ to
$\Gamma(X,\zeta)$, where $\Gamma(X,\zeta): L^2(P)\times\Z\to L^2(P)$
is a given function. In many applications, the process shifts as a result of a simple translation, i.e. $\Gamma(y,\zeta)=y+\zeta$.  \\Simultaneously, the conditional law jumps from $\mu_{t}=\L(X|\F_t^{(1)} )$ to 
\begin{align} \label{muG}
\mu_t^{\Gamma(X,\zeta)}:=\L(\Gamma(X,\zeta)|\F_t^{(1)}).
\end{align}
An {\em impulse control}\index{impulse control} for this system is a double (possibly
finite) sequence
\[
v=(\tau_1,\tau_2,\ldots,\tau_j,\ldots;
\zeta_1,\zeta_2,\ldots,\zeta_j,\ldots)_{j\leq M},\quad
M\leq\infty,
\]
where $0 \leq \tau_1 \leq \tau_2 \leq \cdots$ are $\F_t$-stopping
times (the {\em intervention times})\index{intervention time}
and $\zeta_1,\zeta_2,\ldots$ are the corresponding {\em
impulses}\index{impulses} at these times. Mathematically, we assume that $\tau_j$ is a stopping time with respect to a suitable filtration $\{\F_t\}_{t\ge0}$, with $\tau_{j+1}\ge\tau_j$
 and 
$\zeta_j$ is $\F_{\tau_j}$-measurable for all $j$.
We let $\V$ denote the set of all impulse controls.\\

If $v=(\tau_1,\tau_2,\ldots;\zeta_1,\zeta_2,\ldots) \in \V$, the corresponding state process  $Y^{(v)}(t)$ is
defined by
\begin{align}
& Y^{(v)}(0^-) = y \quad \mbox{and} \quad  Y^{(v)}(t) = Y(t);\quad 0< t\leq \tau_1, \label{513}\\[.5ex]
& Y^{(v)}(\tau_j) = \Big(\tau_j,\Gamma[\check{X}^{(v)}(\tau_j^-),\zeta_j],\L(\Gamma[\check{X}^{(v)}(\tau_j^-),\zeta_j]|\F_t^{1})\Big), \quad j=1,2,\ldots  \label{514}\\[1ex]
&   {\rm d}Y^{(v)}(t) = F(Y^{(v)}(t)){\rm d}t
    +G(Y^{(v)}(t)){\rm d}B(t) \nonumber \\ 
&\qquad\quad\hspace*{17pt} +\int_{\RB^k} H(Y^{(v)}(t^-),z)\widetilde{N}({\rm
d}t,{\rm d}z) \quad \label{515}
     \hbox{for\;$\tau_j< t<\tau_{j+1}\wedge \tau^\ast$}, 
\end{align} 
where
 we have used the notation
     \begin{equation*}
\check{X}^{(v)}(\tau_j^-)=X^{(v)}(\tau_j^-)
+\Delta_N X(\tau_j),
\end{equation*}
$\Delta_N X^{(v)}(t)$ being the jump of $X^{(v)}$ stemming from the jump of the random measure
$N(t,\cdot)$
Note that we distinguish between the (possible) jump of
$X^{(v)}(\tau_j)$ stemming from the random measure $N$, denoted
by $\Delta_N X^{(v)}(\tau_j)$
and the jump caused by the intervention $v$, given by
\begin{equation*}
\Delta_v X^{(v)}(\tau_j):=\Gamma(\check{X}^{(v)}(\tau_j^-),\zeta)
     -\check{X}^{(v)}(\tau_j^-).
\end{equation*}
Accordingly, at the time $t= \tau_j$,  $X^{(v)}(t)$ jumps from
$\check{X}^{(v)}(\tau_j^-)$ to
$\Gamma[\check{X}^{(v)}(\tau_j^-), \zeta_j]$\\
and $\mu_{\tau_j^-}$ jumps to
\begin{align*}
\mu_{\tau_j}=\L(\Gamma[\check{X}^{(v)}(\tau_j^-),\zeta_j]|\F_{\tau_j}^{1}).
\end{align*}
Consider a fixed open set (called the solvency region) $\S\subset [0,\infty) \times \RB^d \times \MM$. It represents the set in which the game takes place since it will end once the controlled process leaves $\S$. In portfolio optimization problems, for instance, the game ends in case of bankruptcy, which may be modelled by choosing $\S$ to be the set of states where the capital is above a certain threshold.
Define
\begin{equation*}
\tau_{\S}=\inf \{ t\in(0,\infty); Y^{(v)}(t)\not\in \S\},
\end{equation*}
and
\begin{equation*}
    \T=\left \{ \tau \; ; \; \textnormal{stopping time,} \; 0 \le \tau \le \tau_\S \right\}.
\end{equation*}

Suppose we are given a continuous {\em profit function}\index{profit rate}
$f:\S\to\RB$ and a continuous {\em bequest function}\index{bequest rate}
$g:\S \to\RB$. Moreover, suppose the profit/utility of
making an intervention with impulse $\zeta\in \Z$ when the state
is $y$ is $K(y,\zeta)$, where $K:\S\times\Z\to\RB$ is a given
continuous function.

We assume we are given a set $\V$ of {\em admissible impulse
controls}\index{admissible impulse controls} which is included in
the set of $v=(\tau_1,\tau_2,\ldots;\zeta_1,\zeta_2,\ldots)$ such
that a unique solution $Y^{(v)}$ of (\ref{513})--(\ref{515})
exist, for all $v\in \V$, and the following additional properties hold, assuring that the performance functional below is well-defined: 
\begin{equation*}
E^y\Big[ \int_0^{\tau_{\S}} f^-(Y^{(v)}(s))ds\Big]
<\infty, \quad \hbox{for all\;$y\in \S$,\;$v\in\V$},
\end{equation*}

\begin{equation*}
E^y\left[ g^-(Y^{(v)}(\tau_{\S}))\mathbbm{1}_{[\tau_{\S}<\infty]}\right]
<\infty, \quad \hbox{for all\;$y\in \S$,\;$v\in\V$},
\end{equation*}
and
\begin{equation*}
E^y\left[ \sum_{\tau_j\leq \tau_{\S}}
K^-(\check{Y}^{(v)}(\tau_j^-),\zeta_j)\right] <\infty, \quad
\hbox{for all\;$y\in\S$,\;$v\in\V$},
\end{equation*}
where ${E}^y$ denotes expectation given that  $Y(0)=y$.\\
We now define the performance criterion, which consists of three parts: a continuous time running profit in $[0,\tau_\S]$, a terminal bequest value  if the game ends, and a discrete-time intervention profit, namely \begin{align*}
J^{(v)}(y)&= E^y \Bigg[ \int_0^{\tau_{\S}} f(Y^{(v)}(t)){\rm d}t
+g(Y^{(v)}(\tau_{\S}))\mathbbm{1}_{[\tau_{\S}<\infty]}+\sum_{\tau_j\leq \tau_{\S}}
K(\check{Y}^{(v)}(\tau_j^-),\zeta_j)\Bigg].
\end{align*}
We consider the following {\em impulse control problem}:

\begin{problem}
Find $\Phi(y)$ and $v^\ast\in\V$ such that
\begin{equation*}
\Phi(y)=\sup \{ J^{(v)}(y);v\in\V\}=J^{(v^\ast)}(y),\quad y \in \S.
\end{equation*}
The function $\Phi(y)$ is called the value function and $v^\ast$ is called an optimal control.
\end{problem}
The following concept is crucial for the solution of this problem.
\begin{definition}\label{def511}
Let $\H$ be the space of all measurable functions $h:\S\to\RB$. 
The {\em intervention operator}\index{intervention operator}
$\M:\H\to\H$ is defined by
\begin{equation}\label{6.1.17}
\M h(s,X,\mu)=\sup_{\zeta \in \Z} \{ h(s,\Gamma(X,\zeta),\mu^{\Gamma(X,\zeta)})+K(y,\zeta), \zeta\in\Z 
\text{ and } (s,\Gamma(X,\zeta),\mu^{\Gamma(X,\zeta)})\in \S\},
\end{equation}
where $\mu^{\Gamma(X,\zeta)}$ is given by \eqref{muG}.
\end{definition}

Let $\mathcal{C}^{(1,2,2)}(\S)$ denote the family of functions $\varphi(s,x,\mu):\S \to \RB$ which are continuously differentiable w.r.t. $s$ and twice continuously Fr\'echet  
differentiable w.r.t. $x \in \RB^d$ and $\mu \in \MM$. We let $\nabla_\mu \varphi \in L(\MM,\RB)$ (the set of bounded linear functionals on $\MM$) denote the Fr\'echet derivative (gradient) of $\varphi$ with respect to $\mu \in \mathbb{M}$. Similarly, $D^2_\mu\varphi$ denotes the double derivative of $\varphi$ with respect to $\mu$ and it belongs to  $ L(\MM \times \MM,\RB)$ (see Appendix for further details). \\The infinitesimal generator $\mathcal{G}$ of the Markov jump diffusion process $Y(t)$ is defined on functions $\varphi \in \mathcal{C}^{(1,2,2)}(\S)$  by  
\begin{align*}
&\mathcal{G} \varphi= \frac{\partial \varphi}{\partial s} +\sum_{j=1}^d \alpha_j \frac{\partial \varphi}{\partial x_j} + \langle \nabla_{\mu} \varphi, A_0^{*} \mu \rangle + \tfrac{1}{2}\sum_{j,n=1}^{d}  (\beta \beta^{T})_{j,n}\frac{\partial ^2 \varphi}{\partial x_j \partial x_n} \nonumber\\
& + \tfrac{1}{2}\sum_{j=1}^d \beta_{j,1}\frac{\partial}{\partial x_j}\langle\nabla_{\mu} \varphi,A_1^{*}\mu\rangle +\tfrac{1}{2} 
\langle A_1^{*}\mu, \langle D_{\mu}^2 \varphi,A_1^{*}\mu\rangle \rangle \nonumber\\
&+\sum_{\ell =1}^k \int_{\RB} \{ \varphi(s, X+\gamma^{(\ell)}, \mu)) - \varphi(s,X,\mu) 
-\sum_{j=1}^d\gamma_j^{(\ell)}  \tfrac{\partial}{\partial x_j} \varphi(s,X,\mu) \}\nu_{\ell}(d\zeta),
\end{align*}
where, as before, 
$A_0^{*}$ is the integro-differential operator
\begin{align*}
A_0^{*}\mu&= -\sum_{j=1}^d D_j [\alpha_j \mu] +\frac{1}{2}\sum_{n,j=1}^d D_{n,j}[(\beta \beta^{(T)})_{n,j} \mu] \nonumber\\
&+\sum_{\ell=1}^k  \int_{\RB}\Big\{\mu^{(\gamma^{(\ell)})}-\mu+\sum_{j=1}^d D_j[\gamma_j^{(\ell)}(s,\cdot,\zeta)\mu]\Big\} \nu_{\ell} \left( d\zeta \right), 
\end{align*}
and
\begin{align*}
A_1^{*}\mu= - \sum_{j=1}^d D_j[\beta_{1,j} \mu].
\end{align*}
We can now state a verification
theorem for conditional McKean-Vlasov impulse control problems,  providing sufficient conditions that a given function is the value function and a given impulse control is optimal. The verification theorem links the impulse control problem to a suitable system of quasi-variational inequalities. \\
Since the process $Y(t)$ is Markovian, we can, with appropriate modifications, use the approach in Chapter 9 in \cite{OS}.\\
For simplicity of notation we will in the following write
\begin{align*}
\overline{\Gamma}(y,\zeta)=(s,\Gamma(x,\zeta),\mu^{\Gamma(x,\zeta)}), \text{ when }  y=(s,x,\mu) \in [0,\infty) \times L^2 (P) \times \MM.
\end{align*}

\begin{theorem}
{Variational inequalities for conditional McKean-Vlasov impulse control}
\label{th512}
\begin{enumerate}
\item[{\rm (a)}] Suppose we can find
$\phi:\bar{\S}\to\RB$ such that
\begin{enumerate}
\item[{\rm (i)}]
$\phi\in \C^1(\S)\cap \C(\bar{\S})$.

\item[{\rm (ii)}]
$\phi\geq\M\phi$ on $\S$.
\\ Define
\[
D=\{y\in \S; \phi(y)>\M\phi(y)\}\quad \hbox{ (the continuation
region\index{continuation region}).}
\]
Assume
\item[{\rm (iii)}]
$\displaystyle \; E^y\left[ \int_0^{\tau_{\S}} Y^{(v)}(t)\mathbbm{1}_{\partial
D}{\rm d}t\right]=0$\quad for all $y\in \S$,\;
$v\in\V$.

\item[{\rm (iv)}]
$\partial D$ is a Lipschitz surface.

\item[{\rm (v)}]
$\phi\in \mathcal{C}^{(1,2,2)}(\S\setminus \partial D)$ with locally bounded
derivatives near $\partial D$.

\item[{\rm (vi)}]
$\G\phi+f\leq 0$ on $\S\setminus \partial {D}$.

\item[{\rm (vii)}]
$\phi(y) = g(y) \text{ for all } y \not\in \S$.

\item[{\rm (viii)}]
$ \{\phi^-(Y^{(v)}(\tau));\tau\in\T\}$ is uniformly
integrable, for all $y\in \S$, $v\in\V$.

\item[{\rm (ix)}]
$\displaystyle E^y\left[|\phi(Y^{(v)}(\tau))| +
\int_0^{\tau_{\S}} |\G \phi(Y^{(v)}(t))| {\rm d}t\right]<\infty$
for all $\tau \in \T, v \in \V, y \in \S$.

\noindent Then
\begin{equation*}
\phi(y)\geq\Phi(y)\quad \hbox{for all $y\in \S$}.
\end{equation*}
\end{enumerate}
\item[{\rm (b)}] Suppose in addition that
\begin{enumerate}
\item[{\rm (x)}]
$\G\phi+f=0$ in $D$.

\item[{\rm (xi)}]
$\hat{\zeta}(y)\in{\rm
Argmax}\{\phi({\overline{\Gamma}(y,\cdot)})+K(y,\cdot)\}\in\Z$ exists for all
$y\in \S$ and $\hat{\zeta}(\cdot)$ is a Borel measurable
    selection.\vskip.2cm

Put $\hat{\tau}_0=0$ and define
$\,\hat{v}=(\hat{\tau}_1,\hat{\tau}_2,\ldots;\hat{\zeta}_1,
\hat{\zeta}_2,\ldots)$ inductively  by

$\hat{\tau}_{j+1}=\inf\{ t>\hat{\tau}_j;
Y^{(\hat{v}_j)}(t)\not\in D\}\wedge \tau_{\S}$ and
$\,\hat{\zeta}_{j+1}=\hat{\zeta}(Y^{(\hat{v}_j)}(\hat{\tau}_{j+1}^-))$

if
$\hat{\tau}_{j+1}<\tau_{\S}$,
where $Y^{(\hat{v}_j)}$ is the result of applying

$\hat{v}_j:=(\hat{\tau}_1,\ldots,\hat{\tau}_j;
\hat{\zeta}_1,\ldots,\hat{\zeta}_j)$ to $Y$.
Suppose
\item[{\rm (xii)}]
$\hat{v}\in\V$ and $\{ \phi(Y^{(\hat{v})}(\tau));
\tau\in\T\}$ is uniformly integrable.

Then
\begin{equation*}
\phi(y)=\Phi(y)\quad \hbox{and $\;\hat{v}\,$ is an optimal
impulse control}.\index{optimal
impulse control}
\end{equation*}
\end{enumerate}\end{enumerate}
\end{theorem}

\noindent

\begin{remark}
    We give the intuitive idea behind intervention operator as in \eqref{6.1.17}: 
    \begin{equation}
\M \Phi(y)=\sup_{\zeta \in \Z} \{ \Phi(\overline{\Gamma}(y,\zeta))+K(y,\zeta), \; \zeta\in\Z 
\text{ and } \overline{\Gamma}(y,\zeta)\in \S\},
\end{equation}
Assume that the value function $\Phi$ is known.
If $y=(s,x,\mu)$ is the current state of the process, and the agent intervenes with impulse of size $\zeta$, the resulting value can be represented as $\Phi(\overline{\Gamma}(y,\zeta))+K(y,\zeta)$, consisting of the sum of the value of $\Phi$ in the new state $\overline{\Gamma}(y,\zeta)$ and the intervention cost $K$. Therefore, $\M\Phi(y)$
represents the optimal new value if the agent decides to make an intervention at $y$.\\
Note that by $(ii)$ $\Phi\geq\M\Phi$ on $\S$, so it is not always optimal to intervene. At the time $\hat{\tau}_j$, the operator should intervene with impulse $\hat{\zeta}_j$ when the controlled process leaves the continuation region, that is when $\Phi(Y^{\hat{v}_j}) \leq \M\Phi(Y^{\hat{v}_j})$. 
\end{remark} 
\dproof 
{(a)} 
By an approximation argument (see e.g. Theorem 3.1 in \cite{OS}) and (iii)--(v), we may assume that
$\phi\!\in\! C^2(\S) \cap C(\bar{\S})$. Choose
$\,v\!=\!(\tau_1,\tau_2,\ldots;\zeta_1,\zeta_2,\ldots)\!\in\!\V$
and set $\tau_0 = 0$. By \nobreak another approximation argument we may
assume that we can apply the Dynkin formula to the stopping times
$\tau_j$. Then for $j=0,1,2,\ldots$, with $Y=Y^{(v)}$
\begin{equation*}
E^y[\phi(Y(\tau_j))]-E^y [\phi(\check{Y}(\tau_{j+1}^-))] =-E^y
\left[ \int_{\tau_j}^{\tau_{j+1}} \G\phi(Y(t)){\rm d}t\right],
\end{equation*}
where $\check{Y}(\tau_{j+1}^-)=Y(\tau_{j+1}^-)+\Delta_N Y(\tau_{j+1})$,
as before. Summing this from $j=0$ to $j=m$ we get
\begin{align}
\phi(y)&+  \sum_{j=1}^m E^y [\phi(Y(\tau_j))
-\phi(\check{Y}(\tau_j^-))]
-E^y[\phi(\check{Y}(\tau_{m+1}^-))] \nonumber \\
&=-E^y \left[ \int_0^{\tau_{m+1}} \G\phi(Y(t)){\rm d}t\right] \geq
E^y \left[ \int_0^{\tau_{m+1}} f(Y(t)){\rm d}t\right].
\label{5123}
\end{align}

\noindent Now
\begin{align*}
\phi(Y (\tau_j))& =\phi(\Gamma(\check{Y}(\tau_j^-),\zeta_j))
\nonumber \\
&\leq\M\phi(\check{Y}(\tau_j^-))-K(\check{Y}(\tau_j^-),\zeta_j)\quad
   \hbox{if $\;\tau_j<\tau_{\S}$\ by\ \eqref{6.1.17}} \end{align*}
\noindent and
   \begin{align*}
  \phi(Y
(\tau_j))&=\phi(\check{Y}(\tau_j^-))\quad
   \hbox{if $\;\tau_j=\tau_{\S}$ \ by (vii).}
\end{align*}
Therefore
\[
\M\phi(\check{Y}(\tau_j^-))-\phi(\check{Y}(\tau_j^-))
\geq \phi(Y(\tau_j))-\phi(\check{Y}(\tau_j^-))
+K(\check{Y}(\tau_j^-),\zeta_j),
\]
and
\begin{align*}
&\phi (y) +  \sum_{j=1}^m E^y
   [\{\M\phi(\check{Y}(\tau_j^-))-\phi(\check{Y}(\tau_j^-)) \}
    \mathbbm{1}_{[\tau_j<\tau_{\S}]}] \\
&\qquad \geq E^y \left[ \int_0^{\tau_{m+1}}
     f(Y(t)){\rm d}t+\phi(\check{Y}(\tau_{m+1}^-))
      +\sum_{j=1}^m K(\check{Y}(\tau_j^-),\zeta_j) \right].
\end{align*}
Letting $m\to M$ and using quasi-left continuity of $Y(\cdot)$, we get
\begin{equation}
\phi(y) \geq E^y\! \left[ \!\int_0^{\tau_{\S}}\! f(Y(t)){\rm
d}t+g(Y(\tau_{\S})) \mathbbm{1}_{[\tau_{\S}<\infty]}
+ \sum_{j=1}^M
K(\check{Y}(\tau_j^-),\zeta_j)\!\right]\!=\!J^{(v)}(y).\label{5125}
\end{equation}
Hence $\phi(y)\geq\Phi(y)$.\\

\noindent {(b)}
Next assume (x)--(xii) also hold. Apply the
above argument to
$\hat{v}=(\hat{\tau}_1,\hat{\tau}_2,\ldots;\hat{\zeta}_1,\hat{\zeta}_2,\ldots)$.
Then by (x) we get {\em equality} in (\ref{5123}) and by our
choice of $\zeta_j=\hat{\zeta}_j$ we have {\em equality} in
(\ref{5125}). Hence
\[
\phi(y)=J^{(\hat{v})}(y),
\]
which combined with (a) completes the proof.

\fproof


\section{Example: Optimal stream of dividends 
under transaction costs} 
In this Section, we solve explicitly an optimal stream of dividends 
under transaction costs. \\ 
To this end, for $v = (\tau_1, \tau_2, \ldots  ;  \zeta_1, \zeta_2, \ldots)$
with $\zeta_i \in \RB_+$, we define \\$Y^{(v)}(t)=(s+t, X^{(v)}(t),\mu_t^{(v)})$ by
\begin{align}
&dX(t)= {E}\left[ X(t)\mid \mathcal{F}_{t}^{(1)}\right] \Big(\alpha
_{0}dt+\sigma _{1}dB_{1}(t)+\sigma _{2}dB_{2}(t) +\int_{\RB}\gamma _{0}(\zeta )\widetilde{N}(dt,d\zeta
)\Big), \label{Xt}\\ \nonumber
&\mu_t^{(v)} = \L(X^{(v)}(t)|\F_t^{(1)}); \quad \tau_i < t < \tau_{i+1},\\  \nonumber
&X^{(v)}(\tau_{i+1})= \check{X}^{(v)}(\tau^-_{i+1}) - (1 +
\lambda) \zeta_{i+1} - c,\\ \nonumber
&\mu_{\tau_{i+1}}^{(v)}= \L(X^{(v)}(\tau_{i+1})|\F^{(1)}_{\tau_{i+1}});
\quad i = 0,1,2, \ldots, \\ \nonumber
&X^{(v)}(0^-)= x > 0; \text{ a.s., }
\end{align}
where $\alpha_0, \sigma_1 \neq 0, \sigma_2 \neq 0, \lambda \geq 0$, and $c >
0$ are constants with $ -1 \leq \gamma_0(z) $
a.s. $\nu$. \\
Here $X(t)$ represents the amount available at time $t$ of a cash flow. We assume that it satisfies the McKean-Vlasov equation in \eqref{Xt}. Note that at any time $\tau_{i}, \; i=0,1,2,\ldots,$ the system jumps from
$\check{X}^{(v)}(\tau_i^-)$ to
$$X^{(v)}(\tau_{i})= \Gamma[\check{X}^{(v)}(\tau_i^-), \zeta_i] = \check{X}^{(v)}(\tau^-_{i}) - (1 +
\lambda) \zeta_{i} - c, $$
 where the quantity $c+\lambda\zeta_i$ represents the transaction cost with a \textit{fixed} part $c$ and a  \textit{proportional} part $\lambda\zeta_i$, while $\zeta_i$ is the amount we decide to take out at time $\tau_i$. \\
 At the same time  $\mu_{\tau_i^-}$ jumps to
\begin{align*}
\mu_{\tau_i}=\L(\check{X}^{(v)}(\tau_i^-)|\F_{\tau_i}^{1}).
\end{align*}

\begin{problem}\label{prob}
We want to find $\Phi$ and $v^\ast \in \V$ such that
\begin{align}\nonumber
\Phi(s,x,\mu) = \sup_v J^{(v)}(s,x,\mu) = J^{(v^\ast)}(s,x,\mu),
\end{align}
where
\begin{align}\nonumber
J^{(v)}(s,x,\mu) = J^{(v)}(y) =E^y \left[\sum_{\tau_k < \tau_{\S}} {\rm e}^{-\rho
(s + \tau_k)} \zeta_k \right] \qquad (\rho > 0 \text{ constant})
\end{align}
is the expected discounted total dividend up to time $\tau_{\S}$, where
\begin{align}\nonumber
\tau_{\S} = \tau_{\S}(\omega) = \inf \{t > 0  ; P^y[ E^y[X^{(v)}(t)|\F_t^{(1)}] \leq
0 ] > 0\}
\end{align}
is the time of bankruptcy. 
\end{problem} 
To put this problem into the context above,  we define
\begin{align*}
&Y^{(v)}(t)=\begin{bmatrix}s+t\\ X^{(v)}(t)\\ \mu^{(v)}_t\end{bmatrix}, \quad  Y^{(v)}(0^-)=\begin{bmatrix} s\\ x\\ \mu\end{bmatrix}=y,\\
&\Gamma(y,\zeta)=\Gamma(s,x,\mu)=(s,x-c-(1+\lambda)\zeta, \L(x-c-(1+\lambda)\zeta)|\F^{(1)}), \quad x \in L^2(P),\\
&K(y,\zeta)={\rm e}^{-\rho s}\zeta, \\
&f\equiv g\equiv 0, \\
&\S =\left\{(s,x,\mu):x>0 \text{ a.s. }\right\}.
\end{align*}
Comparing with our Theorem, we see that in this case we have $d=1,m=2,k=1$ and%
\begin{equation*}
\alpha _{1}=\alpha _{0}\left\langle \mu ,q\right\rangle ,\beta _{1}=\sigma
_{1}\left\langle \mu ,q\right\rangle ,\beta _{2}=\sigma _{2}\left\langle \mu
,q\right\rangle ,\gamma (s,x,\mu ,\zeta )=\gamma _{0}(t,\zeta )\left\langle
\mu ,q\right\rangle ,
\end{equation*}%
where we have put $q(x)=x$ so that $\left\langle \mu_t ,q\right\rangle =E\left[ X(t)\mid 
\mathcal{F}_{t}^{(1)}\right] .$ \\ Therefore the operator $\mathcal{G}$ takes the form
\begin{eqnarray}
\mathcal{G}{\varphi }(s,x,\mu ) &=&\frac{\partial \varphi }{\partial s}+\alpha
_{0}\left\langle \mu ,q\right\rangle \frac{\partial \varphi }{\partial x}%
+\left\langle \nabla _{\mu }\varphi ,A_{0}^{\ast }\mu \right\rangle 
 \nonumber \\
&&+\tfrac{1}{2}(\sigma _{1}^{2}+\sigma _{2}^{2})\left\langle \mu
,q\right\rangle ^{2}\frac{\partial ^{2}\varphi }{\partial x^{2}}+\frac{1}{2}%
\sigma _{1}\left\langle \mu ,q\right\rangle \frac{\partial }{\partial x}%
\left\langle \nabla _{\mu }\varphi ,A_{1}^{\ast }\mu \right\rangle  
\nonumber \\
&&+\tfrac{1}{2}\left\langle A_{1}^{\ast }\mu ,\left\langle D_{\mu
}^{2}\varphi ,A_{1}^{\ast }\mu \right\rangle \right\rangle   \nonumber \\
&&+\int_{\RB}\left\{ \varphi (s,x+\gamma _{0}\left\langle \mu
,q\right\rangle ,\mu )-\varphi (s,x,\mu )-\gamma _{0}\left\langle \mu
,q\right\rangle \frac{\partial }{\partial x}\varphi (s,x,\mu )\right\} \nu
(d\zeta ),  \nonumber
\end{eqnarray}%
where%
\begin{equation*}
A_{0}^{\ast }\mu =-D[\alpha _{0}\left\langle \mu ,q\right\rangle \mu ]+\tfrac{%
1}{2}D^{2}[(\sigma _{1}^{2}+\sigma _{2}^{2})\left\langle \mu ,q\right\rangle
^{2}\mu ],  
\end{equation*}%
and%
\begin{equation*}
A_{1}^{\ast }\mu =-D[\sigma _{1}\left\langle \mu ,q\right\rangle \mu ].
\end{equation*}%
The adjoints of the last two operators are%
\begin{equation*}
A_{0}\mu =\alpha _{0}\left\langle \mu ,q\right\rangle D\mu +\tfrac{1}{2}%
(\sigma _{1}^{2}+\sigma _{2}^{2})\left\langle \mu ,q\right\rangle
^{2}D^{2}\mu,  
\end{equation*}
and%
\begin{equation*}
A_{1}\mu =\sigma _{1}\left\langle \mu ,q\right\rangle D\mu .  
\end{equation*}
In this case the intervention operator gets the form
\begin{equation*}
\M h(s,x,\mu)=\sup \left\{ h(s,x-c-(1+\lambda)\zeta,\mu^{x-c-(1+\lambda)\zeta})
+e^{-\rho t}\zeta; \quad 0\leq\zeta\leq\frac{x-c}{1+\lambda}\right\}.
\end{equation*}
Note that the condition on $\zeta$ is due to the fact that the impulse
must be positive and $x-c-(1+\lambda)\zeta$ must belong to $\S$. We
distinguish between two cases:\\
1. $\alpha_0>\rho$. In this case, suppose we wait until some time $t_1$ and then take out
\[
\zeta_1=\frac{X(t_1)-c}{1+\lambda}.
\]
Noting that $E^y|X(t)]= x \exp(\alpha_0 t)$ for $t<t_1$, we see that the corresponding performance  is
\begin{align*}
J^{(v_1)} (s,x,\mu)&= E^y \Bigg[
\frac{{\rm e}^{-\rho(t_1+s)}}{1+\lambda} (X(t_1)-c)\Bigg] \\[4pt]
&=E^x \Bigg[ \frac{1}{1+\lambda} \big( x{\rm e}^{-\rho s}
{\rm e}^{(\alpha_0-\rho)t_1} - c \text{ } {\rm e}^{-\rho(s+t_1)}\big) \Bigg] \\[4pt]
&\rightarrow \infty\ \text{as}\ t_1 \rightarrow \infty.
\end{align*}
Therefore we obtain $\Phi(s,x,\mu)=+\infty$ in this case.\\
2.  $\alpha_0<\rho$. We look for a solution by using the results of
Theorem~\ref{th512}. \\
We guess that the
continuation region is of the form
\[
D=\left\{ (s,x,\mu): 0< \l \mu,q \r<\bar{x}\right\}
\]
for some $\bar{x} > 0$ (to be determined),
and in $D$ we try a value function of the form
\begin{align*}
\varphi(s,x,\mu)=e^{-\rho s}\psi(\l \mu,q \r).
\end{align*}
This gives

$\G\phi(s,x,\mu)= e^{-\rho s} \G_0\psi(\l \mu,q \r)$, where
\begin{align*}
\G_0 \psi(x,\mu)&= -\rho \psi(\l \mu,q \r) 
+\left\langle \nabla _{\mu }\psi,A_{0}^{\ast }\mu \right\rangle 
 +\frac{1}{2}%
\sigma _{1}\left\langle \mu ,q\right\rangle \frac{\partial }{\partial x}%
\left\langle \nabla _{\mu }\psi ,A_{1}^{\ast }\mu \right\rangle  
\nonumber \\
&+\frac{1}{2}\left\langle A_{1}^{\ast }\mu ,\left\langle D_{\mu
}^{2}\psi ,A_{1}^{\ast }\mu \right\rangle \right\rangle   \nonumber \\
&+\int_{\RB}\left\{ \psi (x+\gamma _{0}\left\langle \mu
,q\right\rangle ,\mu )-\psi (x,\mu )-\gamma _{0}\left\langle \mu
,q\right\rangle \frac{\partial }{\partial x}\psi (x,\mu )\right\} \nu
(d\zeta ). 
\end{align*}
By the chain rule for Fr\'echet derivatives (see Appendix), we have
\begin{align*}
\nabla_{\mu}\psi (h)= \psi'(\l \mu,q \r) \l h,q\r \text{ and } D^2_{\mu} \psi(h,k)= \psi''(\l \mu,q \r)\l h,q\r \l k,q\r.
\end{align*}
Therefore,
\begin{align*}
\l \nabla_{\mu} \psi, A_0^{*}\mu\r = \psi'(\l\mu,q\r)\l A_0^{*}\mu,q\r =\psi'(\l \mu,q\r )\l \mu,A_0q\r = \psi'(\l \mu,q\r )\alpha_0 \l \mu,q\r,
\end{align*}
and similarly
\begin{align*} 
\tfrac{1}{2} \l A_1^{*}\mu,\l D^2_{\mu} \psi,A_1^{*} \mu\r \r &= \tfrac{1}{2} \psi''(\l \mu,q \r \l A_1^{*} \mu,q\r \l A_1^{*} \mu,q \r= \tfrac{1}{2} \psi''(\l \mu,q \r )\l \mu,A_1q\r \l \mu,A_1q \r \nonumber\\
&=\tfrac{1}{2} \psi''(\l \mu,q \r ) \sigma_1^2 \l \mu,q\r^2.
\end{align*}
Moreover, since $\psi $ does not depend on $x$ we see that%
\begin{equation*}
\int_{\RB}\left\{ \phi (s,x+\gamma _{0}\left\langle \mu
,q\right\rangle ,\mu )-\phi (s,x,\mu )-\gamma _{0}\left\langle \mu
,q\right\rangle \frac{\partial \phi }{\partial x}(s,x,\mu )\right\} \nu
(d\zeta )=0. 
\end{equation*}
Substituting this into the expression for $\G_0 \psi$ we get, with $u=\l \mu,q \r$,
\begin{align*}
\G_0 \psi(u)= -\rho \psi (u) + \alpha_0 u \psi'(u) + \tfrac{1}{2} \sigma_1^2 u^2 \psi''(u).
\end{align*}
By condition (x) we are required to have $\G_0 \psi(u)=0$ for all $u \in (0,\bar{x})$,
and this equation has the general solution
\begin{align*}
\psi(u)=C_1 u^{\gamma_1} + C_2 u^{\gamma_2},\quad u \in (0,\bar{x}),
\end{align*}
where $\gamma_1 > 1, \gamma_2 < 0, \text{ and } C_1, C_2$ are constants.
Since we expect $\phi$ to be bounded near 0, we guess that $C_2=0$.\\
We guess that it is optimal to wait till $u=\l\mu_t,q \r= E^y[X(t)|\F_t^{(1)}]$
reaches or exceeds a value $u=\bar{u} > c$ and then take out as
much as possible, i.e., reduce $E^y[X(t)|\F_t^{(1)}]$ to 0. Taking the transaction
costs into account this means that we should take out
$$\hat{\zeta}(u) = \frac{u-c}{1+\lambda} \text{ for } u \geq \bar{u}.$$
We therefore propose that $\psi(u)$ has the form
$$\psi(u) = 
\begin{cases} C_1 u^{\gamma_1} \text{for}\ 0 < u < \bar{u} \\ 
\displaystyle \frac{u-c}{1+\lambda} \text{ for } u \geq \bar{u}.
\end{cases}
$$
Continuity and differentiability of $\psi(u)$ at $u = \bar{u}$
give the equations
$$C_1 \bar{u}^{\gamma_1} = \frac{\bar{u} - c}{1 + \lambda},$$
and
$$C_1 \gamma_1 \bar{u}^{\gamma_1 - 1} = \frac{1}{1 + \lambda}.$$
Combining these we get
$$\bar{u} = \frac{\gamma_1 c}{\gamma_1 - 1}\quad \text{and}\quad C_1 = \frac{\bar{u} - c}{1 + \lambda} \bar{u}^{- \gamma_1}.$$
With these values of $\bar{u}$ and $C_1$, we have to verify that
$\psi$ satisfies all the requirements of Theorem \ref{th512}. We
check some of them:\vspace*{4pt}

\noindent(ii)
 $\phi \geq \M \phi$ on $\S$.\\
In our case we have $\Gamma(s,X,\mu)=(s,X-c-(1+\lambda)\zeta, \mu^{X-c-(1+\lambda)\zeta})$ and hence we get
\begin{align*}
 \M \phi (s,X,\mu)&=  \sup_{\zeta} \Big \{ \phi(s,X - c - (1 + \lambda) \zeta), \mu^{X - c - (1 + \lambda) \zeta} )+e^{-\rho s} \zeta; \ 0 \leq \zeta \leq \displaystyle \frac{\bar{u}-c}{1+\lambda} \Big \}\\
&= e^{-\rho s} \sup_{\zeta} \Big \{C_1 \l \mu^{X-c-(1+\lambda)\zeta},q\r^{\gamma_1} + \zeta; \ 0 \leq \zeta \leq \displaystyle \frac{\bar{u}-c}{1+\lambda} \Big\}\\
&= e^{-\rho s} \sup_{\zeta} \Big\{C_1 (\l \mu,q(x) - c - (1+\lambda)\zeta\r^{\gamma_1} + \zeta; \ 0 \leq \zeta \leq \displaystyle \frac{\bar{u}-c}{1+\lambda}\Big \}\\
&= e^{-\rho s} \sup_{\zeta} \Big\{ C_1 (\l \mu,q \r - c - (1+\lambda)\zeta)^{\gamma_1} + \zeta; \ 0 \leq \zeta \leq \displaystyle \frac{\bar{u}-c}{1+\lambda}\Big\}.
\end{align*}
If $u- c - (1 + \lambda) \zeta \geq \bar{u}$, then
$$\psi (u-c-(1+\lambda)\zeta) + \zeta = \frac{u - 2c}{1 + \lambda} < \frac{u-c}{1 + \lambda} = \psi(u)$$
and if $u-c-(1+\lambda)\zeta < \bar{u}$ then
$$h(\zeta) := \psi(u-c-(1+\lambda)\zeta) + \zeta = C_1 (u-c-(1+\lambda)\zeta)^{\gamma_1} + \zeta.$$
Since
$$h'\left( \frac{u-c}{1+\lambda}\right) = 1\ \text{and}\ h''(\zeta) > 0,$$
we see that the maximum value of $\displaystyle h(\zeta)  ; \ 0 \leq \zeta \leq
\frac{u-c}{1+\lambda}$, is attained at $\displaystyle \zeta = \hat{\zeta}(u) =
\frac{u-c}{1+\lambda}$.

\noindent Therefore
$$\M \psi(u) = \max \left( \frac{x-2c}{1+\lambda},\frac{u-c}{1+\lambda}\right)
 = \frac{u-c}{1+\lambda}\ \text{for all}\ u > c.$$
Hence $\M \psi(u) = \psi(u)$ for $u \geq \bar{u}$.\\
For $0 < u < \bar{u}$ consider
$$k(u) := C_1u^{\gamma_1} - \frac{u-c}{1+\lambda}.$$

\noindent Since
$$k(\bar{u}) = k'(\bar{u}) = 0\quad \text{and}\quad k''(u) > 0\ \text{for all}\ u,$$
we conclude that
$$k(u) > 0\quad \text{for}\ 0 < u < \bar{u}.$$
Hence
$$\psi(u) > \M \psi(u)\quad \text{for}\ 0 < u < \bar{u}.$$

\noindent(vi)
 $A_0 \psi(u) \leq 0$ for $u \in \S \backslash \bar{D}$ i.e., for
$u > \bar{u}$.
For $u > \bar{u}$, we have
{\begin{align*}
A_0 \psi(u) &= - \rho \frac{u-c}{1+\lambda} + \alpha_0 u \text{ }
\frac{1}{1+\lambda}\\
&\quad +\int_{u \kern-0.5pt +\kern-0.5pt  \gamma u z < \bar{u}}\left\{C_1(u + \gamma u z)^{\gamma_1}\kern-0.5pt -\kern-0.5pt \frac{u + \gamma u z - c}{1 + \lambda}\right\}\! \nu({\rm d}z) \\
&\leq  (1+\lambda)^{-1} [(\mu - \rho) u + (\rho + \|\nu\|) c].
\end{align*}}

\noindent Therefore we see that
\begin{align*}
A_0 \psi(u) & \leq 0\ \text{for all}\ u > \bar{u} \\[3pt]
 & \Leftrightarrow (\alpha_0 - \rho) u + (\rho + \|\nu\|) c \leq 0 \text{ for all } u > \bar{u} \\[3pt]
& \Leftrightarrow (\alpha_0 - \rho) \bar{u} + (\rho + \|\nu\|) c \leq 0 \\[3pt]
& \Leftrightarrow \bar{u} \geq \frac{(\rho + \|\nu\|) c}{\rho - \alpha_0} \\[3pt]
& \Leftrightarrow \frac{\gamma_1 c}{\gamma_1 - 1} \geq \frac{(\rho + \|\nu\|) c}{\rho - \alpha_0} \\[3pt]
& \Leftrightarrow \gamma_1 \leq \frac{\rho + \|\nu\|}{\alpha_0 + \|\nu\|}.
\end{align*}
Since
$$F \left(\frac{\rho}{\mu}\right) \geq - \rho + \alpha_0
 \frac{\rho}{\alpha_0} + \frac{1}{2} \sigma^2 \frac{\rho}{\alpha_0} \left(
\frac{\rho}{\alpha_0} - 1 \right) > 0,$$ and $F(\gamma_1) = 0$,
$\gamma_1 > 1$ we conclude that $\gamma_1 < \frac{\rho}{\alpha_0}$ and
hence (vi) holds if $\|\nu\|$ is small enough, say $\|\nu\| \leq K$.\\
Therefore, 
we have 
the following.
\begin{theorem}
Suppose $\| \nu \| \leq K$. Then the value function for \textbf{Problem \ref{prob}} is
$$\Phi(s,x,\mu) = 
\begin{cases} e^{-\rho s} C_1 u^{\gamma_1} \text{for}\ 0 < u < \bar{u} \\ 
\displaystyle e^{-\rho s} \frac{u-c}{1+\lambda} \text{ for } u \geq \bar{u}.
\end{cases}
$$
where $u=\l \mu,q \r=E[X(t) | \F_t^{(1)}]$ and 
$$\bar{u} = \frac{\gamma_1 c}{\gamma_1 - 1}\quad \text{and}\quad C_1 = \frac{\bar{u} - c}{1 + \lambda} \bar{u}^{- \gamma_1}.$$
and $\gamma=\gamma_1 >1$ is the positive solution of the equation
\begin{align*}
 -\rho + \alpha_0 \gamma + \tfrac{1}{2} \sigma_1^2 \gamma (\gamma -1) =0.
 \end{align*}
  The optimal impulse control is to do nothing while $u=E[X(t) | \F_t^{(1)}] < \bar{u}$
  and take out immediately 
  $$\hat{\zeta}(u) = \frac{u-c}{1+\lambda} \text{ when } u \geq \bar{u}.$$
  This brings $E[X(t) | \F_t^{(1)}]$ down to 0, and the system stops. Hence the optimal impulse consists of at most one intervention.
 \end{theorem}

\section{Appendix: Double Fr\' echet derivatives}
In this section we recall some basic facts we are using about the Fr\' echet derivatives of a function $f: V \mapsto W$, where $V,W$ are given Banach spaces.
\begin{definition}
We say that $f$ has a Fr\' echet derivative $\nabla_xf=Df(x)$ at $x \in V$ if there exists a bounded linear map $A:V \mapsto W$ such that
\begin{align*}
\lim_{h \rightarrow 0} \frac{||f(x+h)-f(x)-A(h)||_{W}}{||h||_{V}} =0.
\end{align*}
Then we call $A$ the Fr\' echet derivative of $f$ at x and we put $Df(x) =A$.
\end{definition}
Note that $Df(x) \in L(V,W)$ (the space of bounded linear functions from $V$ to $W$),  for each $x$.

\begin{definition}
We say that $f$ has a double Fr\' echet derivative $D^2 f(x)$ at $x$ if there exists a bounded bilinear map $A(h,k): V \times V \mapsto W$ such that
\begin{align*}
\lim_{k \rightarrow 0} \frac{||Df(x+k)(h)-Df(x)(h)-A(h,k)||_{W}}{||h||_{V}} =0.
\end{align*}
\end{definition}

\begin{example}
\begin{itemize}
\item
Suppose $f:\mathbb{M} \mapsto \RB$ is given by
\begin{align*}
f(\mu)=\<\mu,q\>^2, \text{ where } q(x)=x.
\end{align*}
Then 
\begin{align*}
f(\mu +h) -f(\mu)&= \<\mu+h,q\>^2 -\<\mu,q\>^2\nonumber\\
&= 2 \<\mu,q\> \<h,q\>+\<h,q\>^2, 
\end{align*}
so we see that
\begin{equation*}
Df(\mu)(h)=2\<\mu,q\>\<h,q\>.
\end{equation*}
To find the double derivative we consider
\begin{align*}
&Df(\mu+k)(h)-Df(\mu)(h)\nonumber\\
&=2\<\mu+k,q\>\<h,q\>-2\<\mu,q\>\<h,q\>\nonumber\\
&=2\<k,q\>\<h,q\>,
\end{align*}
and we conclude that
\begin{equation*}
D^2f(\mu)(h,k)=2\<k,q\>\<h,q\>.
\end{equation*}
\item
Next assume that $g:\MM \mapsto \RB$ is given by $g(\mu)=\<\mu,q\>$.
Then, proceeding as above we find that
\begin{align*}
Dg(\mu)(h)&=\<h,q\> \text{ (independent of } \mu)\\
&\text{ and }\nonumber\\
D^2g(\mu) &=0.
\end{align*}

\end{itemize}
\end{example}

\end{document}